\documentclass[10pt]{amsart}
\usepackage{amssymb,amsmath,amsthm,amsfonts}
\usepackage{mathrsfs}
\makeatletter
\def\@fnsymbol#1{\ensuremath{\ifcase#1\or 1\or 2\fi}}
\makeatother

\def\shd{\mathcal{D}}
\def\she{\mathcal{E}}
\def\shm{\mathcal{M}}
\def\sho{\mathcal{O}}
\def\shh{\mathcal{H}}

\newcommand{\C}{\mathbb{C}}

\newcommand{\R}{\mathbb{R}}

\newcommand{\Z}{\mathbb{Z}}
\newtheorem{theorem}{Theorem}[section]
\newtheorem{proposition}[theorem]{Proposition}
\newtheorem{lemma}[theorem]{Lemma}
\newtheorem{corollary}[theorem]{Corollary}

\theoremstyle{definition}
\newtheorem{definition}[theorem]{Definition}

\newtheorem{example}[theorem]{Example}

\newtheorem{remark}[theorem]{Remark}

\begin{document}
\author{Teresa Monteiro Fernandes }

\title[Microsupport of tempered solutions...]{Microsupport of tempered solutions of $\shd$-modules associated to  smooth morphisms}

\date{}
 \begin{abstract}
Let $f:X\to Y$ be a smooth morphism of complex analytic manifolds and let $F$ be an $\mathbb{R}$-constructible complex on $Y$. Let $\shm$ be a coherent $\shd_X$-module. We prove that the microsupport of the solution complex of $\shm$ in the tempered holomorphic functions $t \shh \text{om} (f^{-1}  F, \sho_X)$ is contained in the 1-characteristic variety of $\shm$ associated to $f$, and that the microsupport of the solution complex  in the tempered microfunctions $t\mu hom(f^{-1}F, \sho_X)$ is contained in the 1-microcharacteristic variety of the microlocalization of $\shm$ along $T^*Y\times_Y X$. 
We apply our results to the complex of solutions of  $\shm$ in the sheaf of distributions holomorphic in the fibers of an arbitrary smooth morphism. 
 \end{abstract}

\maketitle

\footnote{The research of the author
was supported by Funda\c c{\~a}o para a Ci\^encia e Tecnologia and Programa Ci\^encia,
Tecnologia e Inova\c c{\~a}o do Quadro Comunit{\'a}rio de Apoio.}
\textit{Dedicated to the memory of my friend and colleague Miguel Ramos, with admiration.}
\section{Introduction}\label{section:intro}
Let $X$ be a complex analytic manifold.  Let $\sho_X$ denote the sheaf of holomorphic functions on $X$ and let $\shd_X$ denote the sheaf of holomorphic differential operators.

  Let $D^b_{\mathbb{R}-c}(\mathbb{C}_X)$ denote the full subcategory of the derived category $D^b(\mathbb{C}_X)$ of the complexes with $\mathbb{R}$-constructible cohomologies (we shall call such a complex an $\mathbb{R}$-constructible sheaf for short).

  Let $t\mathcal{H}\text{om}(.\,\,, \sho_X)$ denote the functor of moderate cohomology introduced  by M. Kashiwara in (\cite{K2} ) and studied by M. Kashiwara and P. Schapira in (\cite{KS3}).  It is a functor from $D^b_{\mathbb{R}-c}(\mathbb{C}_X)$ to $D^b(\shd_X)$ and, in particular, when $M$ is a real manifold and $X$ is a complexification of $M$, by taking as $F$ the dual $D'(\mathbb{C}_M)$ of the constant sheaf on $M$, we get $t\mathcal{H}\text{om}(F, \sho_X) =Db_M$, the sheaf of Schwartz  distributions on $M$.

Also when $F=\C_U$, the constant sheaf on  a subanalytic  Stein open subset $U$,  $t \shh \text{om} (\C_U, \sho_X)$  is concentrated in degree zero (\cite{Be}, Lemma 2.6) and is well understood as the sheaf $\sho^{t-U}_X$ of holomorphic functions which are tempered on $U$.

The notion of microsupport of an object $G$ in $D^b(\mathbb{C}_X)$, $SS(G)$, was
introduced by M. Kashiwara and P. Schapira (\cite{KS1}),  as a subset of the cotangent bundle $T^*X$ and it describes the directions of non propagation for $G$.

 Let $\mathcal{M}$ be a coherent $\shd_X$-module. In other words, locally, $\shm$ is a system of partial differential equations on $X$.
 The estimation of the microsupport of the solutions of $\mathcal{M}$ in  $t\mathcal{H}\text{om}(F, \sho_X)$  (or, more generally, of an object of $D^b_{coh}(\shd_X)$) is a challenging problem in algebraic analysis, since the growth conditions require far beyond algebro-geometrical tools.  One of these tools is the  classical Levi-condition for one operator generalized to systems by the notion of $1$-microcharacteristic variety (of a microdifferential system along an involutive submanifold of $T^*X$) (\cite{L}, \cite{TMF1}, \cite{S}).

 Several particular cases have been settled such  as  the case of a microdifferential operator and the sheaf of tempered microfunctions (also introduced by \cite{A}) in \cite{B}.
   When $\shm$ has real simple characteristics,  the microsupport of the distribution solutions of $\shm$ was studied in \cite{KT} and \cite {AAT}.

 More recently, for a $\C$-constructible object $F$,  we proved in \cite{TMF} that the microsupport of the solutions of  $\shm$ in $t\mathcal{H}om(F, \sho_X)$ is controlled by the 1-microcharacteristic variety of the microlocalised of $\shm$ along any smooth involutive manifold containing the microsupport of $F$. For that we used the crucial fact that, by \cite{K2} (see also\cite{KS3}),  $t \shh \text{om} (F, \sho_X)$ is an object of the bounded derived category of regular holonomic $\shd_X$-modules  which has no counterpart for an arbitrary $\R$-constructible $F$.

In \cite {KMS}, assuming that $\shm$  is regular along a sub-bundle $V$ of $T^*X $  in the sense of \cite{KO}, for any $F\in D^b_{\mathbb{R}-c}(\mathbb{C}_X)$, one obtained the estimate
$$SS(\text{R}\shh\text{om}_{{\shd}_X} (\shm , t \shh \text{om} (F, \sho_X)))\subset V\hat{+} SS(F)^a$$
($a$ denotes  the antipodal map on $T^*X$ and the operation $\hat{+}$ was defined in \cite{KS1}).

To treat the case of a general $\R$-constructible $F$,  
we have tools which allow to reduce to  $F=\C_U$ for a relatively compact subanalytic open set $U$, but $t \shh \text{om} (\C_U, \sho_X)$ is not easy to manipulate and, and, even in the case  $U$ Stein subanalytic, the usual techniques require a  Cauchy theorem which is unavailable for general $U$.

 In this paper we consider the case of $ t \shh \text{om} (  f^{-1}F, \sho_X)$
 where $f:X\to Y$ is a smooth morphism of complex manifolds and $F\in D^b_{\R-c}(\C_Y)$.
To make it clear, realizing $f$ locally as a projection, it means that we deal with growth conditions on the factor $Y$.  This appears to be the farthest one can go  with this kind of techniques in the present state of art.

The control of the microsupport is then expressed in terms of  the 1-characteristic variety (of a coherent $\shd_X$-module with respect to a smooth morphism $f:X\to Y$), introduced in \cite{S}, where it is denoted by $C^1_f(\shm)$. It is a variant of the the notion of $1$-microcharacteristic variety for $\shd_X$-modules and it is a conic closed involutive analytic subset of the relative conormal bundle $T^*(X|Y)$.

 More precisely, denoting by $\tilde{f}$ the canonical projection $T^*X\to T^*(X|Y)$,
  we prove in Theorem \ref{T:1}  that, for any  $G\in D^b_{\mathbb{R}-c}(\mathbb{C}_X)$ such that $SS(G)\subset X\times_Y T^*Y$, we have
 \begin{equation}\label{E:60}
 SS(\text{R}\shh\text{om}_{{\shd}_X} (\shm , t \shh \text{om} (G, \sho_X)))\subset \tilde{f}^{-1}(C^1_f(\shm)).
 \end{equation}
Of course, when $f$ is the identity, the estimate above is trivial   since $$T^*(X|X)=0.$$
Note that,  by  Proposition 5.4.5 of \cite{KS1},  any $G$ satisfying the assumption of Theorem \ref{T:1}  is, locally on $X$, isomorphic to $f^{-1} F$ for some $F\in D^b_{\mathbb{R}-c}(\mathbb{C}_Y)$.

Our main tools for the proofs are:
\begin{itemize}
 \item{The  analysis developped by Jean-Michel Bony in \cite{B}, which we partially adapted in order to obtain a Cauchy-Kowalevskaia theorem of precise type with data satisfying tempered growth conditions on a subanalytic open subset of the form $f^{-1}\Omega$.}
\item{The subanalytic version of Grauert's theorem obtained in \cite{BMF} which asserts that any subanalytic open subset of a real analytic manifold admits, in any of its complexifications, a fundamental system of Stein subanalytic open neighborhoods.}

\end{itemize}
The paper is organized as follows:

In the Section ${2}$ we pass under review the basic material for this work.

In Section ${3}$ we  state our main result Theorem \ref{T:1} and treat a few non trivial examples. In Section ${4}$ we study the realification procedure.

In Section ${5}$ we study the Cauchy problem and propagation in order to prove Theorem \ref{T:1}.
Let us give the guiding ideas of the proof:

We start by treating the Cauchy problem for $ \sho_X^{t-f^{-1}(\Omega)}$, for a given open subanalytic set $\Omega \subset Y$.
More precisely, given  a coherent $\shd_X$-module  $\shm$  and  a closed submanifold non 1-characteristic for   $\shm$  with respect to $f$,
we prove a  Cauchy-Kowalevskaia-Kashiwara theorem (Theorem \ref{T:401})
with data in
 $ \sho_X^{t-f^{-1}\Omega}|_H$ and  $ \sho_H^{t-f^{-1}\Omega\cap H}$.

 This is performed by  reduction of  $f$ to the case of a  projection
  $f: X=Z\times Y\to Y$, and to $\shm=\shd_X/\shd_X P$ for an operator $P$.
Then, given  $H$ a submanifold of $X$ non 1-characteristic for $P$ with respect to $f$, we obtain a precised Cauchy-Kowalevskaia theorem  in Proposition \ref{P:5}. The proof is essentially an adaptation of the Cauchy theorem 2.4.3 in \cite{B}(and of the techniques in \cite{BS})  with the variable in $Y$ playing the role of a parameter.

After obtaining a Zerner type Lemma (cf. Lemma \ref{L:2}), we first prove estimate (1) in the case $G=\mathbb{C}_{\Omega}$,  for $\Omega$   subanalytic open Stein in $Y$, which implies that (1) holds for $G=\mathbb{C}_S$ provided that
\begin{itemize}
  \item{$S$ is a closed subanalytic subset of $Y$ complementary to a finite union of Stein subanalytic open subsets.}

\end{itemize}

 The proof of  Theorem \ref{T:1} is then reduced, by the realification procedure,  to the case $G=f^{-1}\C_U$, for an open subanalytic relatively compact subset $U$ in a real manifold $N$ complexified by $Y$. At this stage, the result follows by the subanalytic version of Grauert's theorem (\cite{BMF}).

As an example of application we obtain (cf. Corollary \ref{T:2}):

Let $X$ be a complex manifold, let $Y$ be a d-dimensional manifold complexifying a real analytic submanifold $N$ and let $f: X\to Y$ be a smooth morphism. Then $t \shh \text{om} (D'(f^{-1}\mathbb{C}_N), \sho_X)$ is concentrated in degree $d$ and
$\shh^d ( t \shh \text{om} (D'(f^{-1}\mathbb{C}_N), \sho_X))$ can be understood as a sheaf of distributions with holomorphic parameters which we shall denote, for the sake of simplicity, $Db_{X|N}$. By Theorem \ref{T:1} we get:
 $$SS(\text{R}\shh\text{om}_{{\shd}_X} (\shm ,Db_{X|N}))\subset \tilde{f}^{-1}(C^1_f(\shm)).$$
 The Appendix is devoted to explain how far we can go with our results in  the microlocal setting. More precisely, although persuaded that the estimation in  Theorem \ref{T:1} holds when we deal with microdifferential systems and the tempered microlocalization functor $t\mu hom$, we were not able to prove it at present. Recall that recently S. Guillermou  proved in \cite{SG} that  $t\mu hom(F,\sho_X)$ is an object of the derived category of $\she_X$-modules.

 Instead, similarly to \cite{KMS}, we can prove it for  the microlocalization of a $\shd_X$-module $\shm$, that is, the $\she_X$-module obtained by tensorizing  $\shm$ by the sheaf $\she_X$ of microdifferential operators.
  More precisely, we obtain  in Theorem \ref{T:300} a variant of the differential estimate with the 1-characteristic variety replaced by  the 1-microcharacteristic variety:

 Let $\she_X$ denote the sheaf of microdifferential operators on $X$ and let $\pi:T^*X\to X$ be the projection.
Let $t\mu hom(., \sho_X)$ denote the tempered microfunctions functor introduced in  \cite{A}, from $D^b_{\R-c}(X)$ to $D^b(\she_X)$.
Setting $$\tilde{\shm}:=\she_X\otimes_{\pi{-1}\shd_X} \pi^{-1}\shm,$$
 noting $V$  the regular involutive submanifold $X\underset{Y}{\times}T^*Y$ of $T^*X$, noting $$\rho_V:T(T^*X)\to T_V(T^*X)$$ the canonical projection, and noting $$C^1_{V^{\cdot}}(\tilde{\shm})$$ the $1$-microcharacteristic variety of $\tilde{\shm}$ along $V^{\cdot}=V\setminus T_X^*X$  (recall that $C^1_{V^{\cdot}}(\tilde{\shm})$ is a subset of $T_V(T^*X)$) where we identify $T^*(T^*X)$ to $T(T^*X)$ by the Hamiltonean morphism, we obtain in Theorem \ref{T:300},
$$SS(R\shh om_{\she_X}(\tilde{\shm}, t\mu hom(G, \sho_X)))\subset \rho_V^{-1}(C^{1}_{V^{\cdot}}(\tilde{\shm})),$$
for any $G\in D^b_{\R-c}(X)$ whose microsupport is contained in $V$.

We thank  M. Kashiwara and P. Schapira for  their useful comments, and also G. Morando, L. Prelli, A.R. Martins, S. Guillermou, for their interest. Last but not the least, I thank J-M Bony for his patience in reading and clearing what I took from his techniques.
\section{Background}
We will mainly follow the notations in ~\cite{KS1}, \cite{KS3}  and in 1.3 of \cite{S}.

On a complex manifold $X$, we consider the sheaf $\mathcal{O}_X$
of holomorphic functions, the sheaf  $\shd b_X$ of Schwartz  distributions on X and the sheaf $\mathcal{D}_X$ of linear
holomorphic differential operators of finite order. For a non negative integer $m$, we note $\shd_X(m)$ the subsheaf of operators of order at most $m$.
For $\shm\in D^b_{coh}(\shd_X)$, we denote by $Char(\shm)$ its characteristic variety.

We denote by $Mod(\shd_X)$
the abelian category of $\shd_X$-modules
and denote by $Mod_{coh}(\shd_X)$
 the abelian category of coherent $\shd_X$-modules.  

Given a morphism $f:X\to Y$ of complex manifolds we shall consider the associated (derived)  functors $\underline{f}^{-1}$ from $D^b(Mod(\shd_Y))$ to $D^b(Mod(\shd_X))$  and $R\underline{f}*$,  from
$D^b(Mod(\shd_X))$ to $D^b(Mod(\shd_Y))$.

\subsection{Microlocal tools and review on microsupport}
Let $X$ be a real analytic manifold. Let
$TX$ be the tangent bundle to $X$ and
$T^*X$ the cotangent bundle, with the projection $\pi:T^*X\to X$. We identify $X$ with
the zero section of $T^*X$.
Given a smooth submanifold $Y$ of $X$, let
$T_YX$ be the normal bundle to $Y$ and let $T_Y^*X$ be the conormal
bundle.


For a morphism $f:X\rightarrow Y$ of manifolds, we denote by
\begin{center}
$f_\pi: X\times_Y
T^*Y\rightarrow T^*Y$ and $f_d:X\times_Y
T^*Y\rightarrow T^*X$
\end{center}
the associated morphisms.
Recall that when $f$ is smooth, $f_d$ is an embedding, and that when $f$ is an embedding, $f_{\pi}$ is smooth.

Let us assume that $f$ is smooth.
Set  $V =X\times_Y T^*Y$. Then  $f_{\pi}$ is smooth and, by $f_d$, $V$ is a sub-bundle of $T^*X$. Recall that the relative cotangent bundle is defined by the exact sequence $$0\to X\underset{Y}{\times}T^*Y\underset{f_d}{\to} T^*X\underset{\tilde{f}}{\to} T^*(X|Y)\to 0.$$
\begin{remark}\label{R:200}
Let $X, Y$ be real manifolds, and let $p: X\times Y\to X$, $q: T^*(X\times Y)\to T^*X$ be the projections. Identifying $Y$ with the zero section of $T^*Y$, we have $T^*X\times Y\subset T^*(X\times Y)$. Let $A\subset T^*X\times Y$. Then $p_{\pi}p_d^{-1}(A)=q(A)$.
\end{remark}

For a subset $A$ of $T^*X$, we denote by $A^a$ the image of $A$ by
the antipodal map $$a:(x,\xi)\mapsto (x;-\xi).$$

The closure of $A$ is denoted by $\overline{A}$. Let $\gamma\subset TX$ be a cone; the polar cone $\gamma^\circ$ to $\gamma$ is
the convex cone in $T^*X$ defined by
\begin{center}
$\gamma^\circ=\{(x;\xi)\in TX; x\in\pi(\gamma)$ and $Re \langle
v,\xi\rangle\geq 0$ for any $(x;v)\in\gamma \}.$
\end{center}

Given $A$ and $B$ two closed $\R^+$-conic subsets of $T^*X$ one associates a closed $\R^+$-conic subset $A\hat{+} B$ of $T^*X$ (cf. Definition 6.2.3 and Proposition 6.2.4 of \cite {KS1}).

Let $j :X\to Y$ be an embedding  of manifolds. The notion of
$j^{\sharp}$ will be used in the Appendix. Recall that it is a  correspondence introduced in \cite{KS1} associating
conic subsets of $T^*X$ to conic subsets of $T^*Y$. It is characterized by the following result:

\begin{lemma}[cf. Remark 6.2.8 of \cite{KS1}]\label{P:132}
Let $\Lambda$ be a conic subset of $T^*Y$.
Let $(x^{\prime}, x^{\prime\prime})$  be a system of local coordinates on
$Y$ such that $X=\{x^{\prime}=0\}$. Let $(x^{\prime},x^{\prime\prime};\xi^{\prime},\xi^{\prime\prime})$  be the
associated coordinates on $T^*Y$. Then

$(x^{\prime\prime}_0; \xi^{\prime\prime}_0)\in j^{\sharp}(\Lambda)$ if and only if there exists a
sequence $\{(x^{\prime}_n,x_n^{\prime\prime};\xi^{\prime}_n,\xi_n^{\prime\prime})\}_n$ in $\Lambda$  such that $$ x^{\prime}_n\xrightarrow[n]{}0,
x^{\prime\prime}_n\xrightarrow[n]{} x^{\prime\prime}_0,\xi_n^{\prime\prime}\xrightarrow
[n]{}\xi_0^{\prime\prime}, \vert x^{\prime}_n\vert \vert {\xi_n}^{\prime}\vert \xrightarrow[n]{} 0. $$
\end{lemma}

 We denote by  $D^b(\mathbb{C}_X)$ the triangulated category  of complexes of $\mathbb{C}$-vector spaces with bounded
cohomologies and by $ D^b_{\mathbb{R}-c}(\mathbb{C}_X)$ (respectively  $ D^b_{\mathbb{C}-c}(\mathbb{C}_X)$) the triangulated category of complexes with bounded and $\mathbb{R}$-constructible cohomologies (respectively with $\mathbb{C}$-constructible cohomologies).

For any $F\in D^b(\mathbb{C}_X)$,  its
microsupport  $SS(F)$,  introduced by M. Kashiwara and P. Schapira in \cite{KS1}, is a closed $\R^+$-conic involutive subset of $T^*X$ given by:
\begin{definition}
Let $F\in D^b(\mathbb{C}_X)$. Let $p\in T^*X$. Then $p\notin SS(F)$ if and only if
there exists an open conic
neighbourhood $U$ of $p$ such that for any $x\in \pi(U)$ and any
$\mathbb{R}$-valued $C^\alpha$-function $\varphi$ defined on a
neighbourhood of $x$ such that $\varphi(x)=0$, $d\varphi(x)\in U$,
one has
\begin{center}
$H^j_{\{\varphi \geq 0\}}(F)_x=0$ for any $j.$
\end{center}

\end{definition}

We shall need need the following description of $SS(\cdot)$:

Given $(x_{0},\xi_{0})\in \mathbb{R}^n\times(\mathbb{R}^n)^*$ and
$\varepsilon\in\mathbb{R}$ we set:
$$H_{\varepsilon}(x_{0},\xi_{0})=\{x\in\mathbb{R}^n;\langle
x-x_{0},\xi_{0} \rangle> -\varepsilon\},$$ and if there is no risk
of confusion we will write $H_{\varepsilon}$ instead of
$H_{\varepsilon}(x_{0},\xi_{0})$.

\begin{proposition}[\cite{KS1}]\label{P:165}
 Let $\Omega$ be a local chart in a neighbourhood of $p$ such that $\Omega$ is identified to an open subset of $\mathbb{R}^n$ and
$p=(x_{0},\xi_{0})$. Assume that $p\notin SS(F)$. Then there exist a proper closed convex cone
$\gamma\subset\mathbb{R}^n$, $\varepsilon>0$ and an open
neighbourhood $\omega$ of $x_{0}$ with $\xi_{0}\in Int(\gamma^\circ)$
such that $(\omega+\gamma^a)\cap \overline{H_{\varepsilon}}\subset X$
and
\begin{center}
$H^j(X; \mathbb{C}_{(x+\gamma^a)\cap H_{\varepsilon}}\otimes
F)=0,$ for any $j\in\mathbb{Z}$, $x\in \omega.$
\end{center}

\end{proposition}

\subsection{Review on $t\shh om$}

 The functor of moderate cohomology from $D^b_{\R-c}(\C_X)$ to $D^b(\shd_X)$,
 $t \shh \text{om} (\cdot, \sho_X)$, was introduced in \cite {K2}. For a detailed study we refer to \cite{KS3}.
 Recall that, when $X$ is a complexification of a real analytic manifold $M$ and $F=D^{\prime}(\C_M)$, we have $t\shh om(F, \sho_X)\simeq \shd b_M$.

 Given an open subanalytic subset $U$ of $X$, one defines the sheaf of tempered holomorphic functions on $U$, $\sho_X ^{t-U}$, by setting:

\noindent $$\Gamma(\Omega; \sho^{t-U})=\{f\in\Gamma(\Omega\cap U; \sho_X); \text{for any compact K}\subset \Omega\,\,
\text{there exist an}$$
$$\text{integer k}\geq0\,\,\text{and a number}\,\,
 C_K\geq0\,\,  \text{such that}\underset{K\cap U}{sup} |f(z)d(z, \complement U)^k|\leq C_K\}.$$

As a byproduct of Section 10.2 of \cite{KS3} we get
\begin{proposition}
 For an open subanalytic subset $U$ of $\mathbb{C}^n$, one has
$$\shh^0 t \shh \text{om} (\mathbb{C}_U, \sho_X)=\sho_X^{t-U}. $$

\end{proposition}

\begin{proposition}[cf \cite{Be}, Lemma 2.6, \cite{Ho}, Theorem 2.5] \label{P:300}

Assume that $U$ is Stein relatively compact. Then
$t \shh \text{om} (\mathbb{C}_U, \sho_X)$ is concentrated in degree 0 and, for any open Stein subset $V$, for any $j\geq 1$,
\begin{equation}\label{E:18}
 H^{j}(V; t\shh  \text{om} (\mathbb{C}_U, \sho_X))=0.
\end{equation}
\end{proposition}

\subsection{Review on the 1-characteristic variety}

Let us now consider a smooth morphism $f:X\to Y$ of complex manifolds, let $\shm\in Mod_{coh}(\shd_X)$ and let us recall the construction of $C^1_f(\shm)$ due to P. Schapira in \cite{S}.
The ring of differential operators relative to $f$ is
$$\shd_{X|Y}:=\{P\in\shd_X; [P,f^{-1}\sho_Y]=0\}$$
which is endowed with the filtration by the order induced by the order on $\shd_X$. We  note $\sigma^1_{X|Y}$ the associated principal symbol regarded as a holomorphic function on $T^*(X|Y)$. We shall keep the notation $\sigma(P)$ for the usual principal symbol of $P$ as a section of $\shd_X$.

The 1-characteristic variety of $\shm$ with respect to $f$ is a closed conic analytic subset of $T^*(X|Y)$ which,
assuming  that $\shm=\shd_X/\mathcal{J}$ for a coherent ideal $\mathcal{J}$ of $\shd_X$, is given by:
$$C^1_f(\shm)=\{\theta\in T^*(X|Y); \sigma^1_{X|Y}(P)(\theta)=0, \forall P\in \mathcal{J}\cap \shd_{X|Y}\}.$$
\begin{proposition}\label{P:151}(\cite{S})
If $0\to \mathcal{L}\to \shm \to \mathcal{N} \to 0$ is an exact sequence of coherent $\shd_X$-modules, then
$$C^1_f(\shm)=C^1_f(\mathcal{L})\cup C^1_f(\mathcal{N}).$$
\end{proposition}

We can extend this construction to $D^b_{coh}(\shd_X)$, setting, for $\shm\in D^b_{coh}(\shd_X)$,  $$C^1_f(\shm):=\cup_j C^1_f(\shh^j(\shm)).$$

\begin{lemma}\label{L:152}
Let $X,Y,Z$ be complex manifolds and let $f:X\to Y$ be a smooth morphism. Let  $p:X\times Z\to X$ denote the projection and  let $q:T^*(X\times Z)\to T^*X$  denote the canonical projection. Then, for any $\shm\in D^b_{coh}(\shd_X)$, one has  $ q^{-1}
\tilde{f}^{-1}(C^1_f(\shm))=(\widetilde{f\times Id_Z})^{-1}(C^{1}_{f\times Id_Z}(\underline{p}^{-1}\shm)).$
\end{lemma}
\begin{proof}
By classic tools, we  may assume that $\shm$ is concentrated in degree zero and is of the form $\shm=\shd_X/\mathcal{J}$, for a coherent ideal $\mathcal{J}$ of $\shd_X$. Let us consider a system of local coordinates $(x)$ on $X$, $(z)$ on $Z$ and $(y)$ on $Y$.
Then $$(\widetilde{f\times Id_Z})^{-1}(C^{1}_{f\times Id_Z}(\underline{p}^{-1}\shm))=\{(\theta, \theta^{\prime})\in T^*X\times T^*Z; \sigma(P)(\theta,\theta^{\prime})=0\}$$  where $P$ varies in $\shd_{X\times Z|Y\times Z}\cap \langle \shd_{X\times Z}\mathcal{J},\shd_{X\times Z} D_z\rangle$ and the result follows by the equality $\shd_{X\times Z|Y\times Z}\cap \langle \shd_{X\times Z}\mathcal{J},\shd_{X\times Z} D_z\rangle=\shd_{X\times Z|Y\times Z}\cap \shd_{X\times Z}\mathcal{J}.$
\end{proof}

For a submanifold $H\subset X$, we shall say that $H$ is non 1-characteristic for $\shm$ with respect to $f$ if $T^*_H X\cap \tilde{f}^{-1}(C_f^1(\shm))\subset T^*_H H.$
\begin{remark}
Let $\she_X$ denote the sheaf of microdifferential operators on $X$. Denoting by $V$ the sub-bundle $X\underset{Y}{\times}T^*Y$ of $T^*X$, as we shall see in Section 6, $C^1_f(\shm)$ is the differential version of the 1-microcharacteristic variety of the microdifferential system $\she_X\underset{\pi^{-1}\shd_X}{\otimes}\pi^{-1}\shm$ along $V\setminus T^*_XX$.
\end{remark}
\section{Main result and examples}
We now state the main result of this work. Its proof will be given in several steps throughout  sections 4 and 5.
\begin{theorem}\label{T:1}
Let $X$ and $Y$ be complex manifolds, let $f: X\to Y$ be a smooth morphism. Then, for any $G\in
D^b_{\mathbb{R}-c}(\mathbb{C}_X)$ such that  $SS(G)\subset X\times_Y T^*Y$ and any $\shm\in Mod_{coh}(\shd_X)$, we have
\begin{equation}\label{E:1}
SS(\text{R}\shh\text{om}_{{\shd}_X} (\shm , t \shh \text{om} (G, \sho_X)))\subset \tilde{f}^{-1}(C^1_f(\shm)).
\end{equation}
\end{theorem}
We illustrate our result with simple examples:
\begin{example}
 Let $X=\C^2=\C_{z}\times\C_{y}$ and let $f: X\to\C_y$ be the projection. Let $P(z,D_z,y)\in \shd_{X|Y}$ be the operator $D_z-y$. Consider the $\shd_X$-module $\shm=\shd_X/\langle P\rangle.$ Then, in the canonical symplectic coordinates $(z,y;\zeta, \eta)$ in $T^*X$, we get $C^1_f(\shm)=\{\zeta=0\}$ hence $(0,0;1,0)\notin C^1_f(\shm)$. Let $F=\C_U$ with $U=\{y; \Re y>0\}$  and let $H_+=\{(z,y);\Re z>0\}$. We have $t\shh om(f^{-1}\C_U, \sho_X)=\sho^{t-\C\times U}$.

Since $e^{-yz}(D_z-y)e^{yz}=D_z$, $\shm$ is regular along $V=\{(z,y;\zeta, \eta); \zeta=0\}$ in the sense of \cite{KO}, hence this example is also an application of \cite{KMS}.
Let $B_{\epsilon}$ be an open polydisc of $X$ with center $(0,0)$ and radius $1\gg\epsilon>0$.
Let $g\in \Gamma(B_{\epsilon}\cap H_+; t\shh om(f^{-1}\C_U, \sho_X))$ (respectively $g\in \Gamma(B_{\epsilon}; t\shh om(f^{-1}\C_U, \sho_X))).$ Then, equation $Pu=g$ has a solution $u$ in $\Gamma(B_{\epsilon}\cap H_+\cap \C\times U; \sho_X)$ (respectively in  $\Gamma(B_{\epsilon}\cap \C\times U;  \sho_X))$ since solving $(D_z-y)f=g$ is equivalent to solve $D_z (e^{-zy}f)=e^{-zy}g$. In both cases it becomes clear that in fact $u$ defines a section of $t\shh om(\C\times U, \sho_X)$.

Suppose now that  $u\in \Gamma(B_{\epsilon}\cap H_+; t\shh om(f^{-1}\C_U, \sho_X))$ is such that $Pu$ extends, for some $\epsilon^{\prime}$, to a section in  $\Gamma(B_{\epsilon^{\prime}}; t\shh om(f^{-1}\C_U, \sho_X))$. Since $u$ extends in the $z$ variable by classical Zerner's Lemma, and  the general solution of the homogeneous equation $Pu=0$ in  $t\shh om(f^{-1}\C_U, \sho_X))$ is of the form $C(y) e^{yz}$, with $C(y)\in t\shh om(\C_U,\sho_Y)$, the extension of $u$  defines a section in $\Gamma(B_{\epsilon^{\prime}}; t\shh om(f^{-1}\C_U, \sho_X))$.
\end{example}
\begin{example}
 In the situation of 1. let $Q(z,D_z,y)=z^2 D_z+y^2$ and let $\mathcal{N}=\shd_X / \langle Q \rangle$. Let $0<\epsilon\ll1$.
Then $C^1_f (\mathcal {N})=\{z=0\}\cup \{\zeta=0\}.$  Therefore $(0,0;1,0)\in \tilde{f}^{-1}(C^1_f(\mathcal{N}))$.  Let $U\subset Y$ be equal to $\{\Re y>0\}$. Consider a solution $u=C(y) e^{-y^2/z}$ of the homogeneous equation, with $C(y)\in \Gamma(B_{\epsilon}, t\shh om(\C_U,\sho_Y))$. Clearly, $u\in \Gamma(\{\Re z>0\}\times B_{\epsilon}(0), t\shh om(f^{-1}\C_U, \sho_X))$. However  $u$ does not extend to any open set containing points of the form $(0,y)$.
\end{example}
\begin{example} Let $X=\C_{z_1}\times\C_{z_2}\times \C_y$, let $P(z_1,z_2,y; D_{z_1}, D_{z_2})=D_{z_1}-yD_{z_2}\in \shd_{X|Y}$ and let $\shm=\shd_X/\langle D_{z_1}-yD_{z_2}\rangle$. Let $f:X\to \C_y$ be the projection and let $U=\{y; \Re y>0\}$ as in the previous examples.
Then $\tilde{f}^{-1}(C^1_f(\shm))=\{(z_1, z_2, y; \zeta_1,\zeta_2,\eta); \zeta_1=y \zeta_2\}$, hence $(0,1,0;1,0,0)$ is 1-non characteristic for $\shm$ with respect to $f$. Given $0<\epsilon\ll1$, the solution  of the Cauchy problem $Pu=0$, $u|_{z_1=0}=1/ye^{1/z_2}$, for $|z_1|<\epsilon$, $\Re y>0, |y|<\epsilon$ and $|z_2-1|<\epsilon$ is $u=1/ye^{\frac{1}{yz_1+z_2}}$ which is a section of $t\shh om(f^{-1} U, \sho_X)$ since $e^{\frac{1}{yz_1+z_2}}$ it is holomorphic for $|z_1|<\epsilon$, $|y|<\epsilon$ and $|z_2-1|<\epsilon$ .
\end{example}

Assume that $Y$ is a complexification of a real analytic manifold $N$ and let $d$=dim $Y$. For a smooth morphism $f:X\to Y$ the complex
$t \shh \text{om} (\mathbb{C}_{f^{-1} (N)},\sho_X)$ is concentrated in degree $d$. In particular, when $f$ is a projection, we can describe
$\shh^d(t \shh \text{om} (\mathbb{C}_{f^{-1} (N)},\sho_X))\otimes or_{f^{-1}(N)|X}$ as a sheaf of distributions with holomorphic parameters which we denote by $Db_{X|N}$ for short. We get the following estimate, which is to relate with the results in \cite{B}.

\begin{corollary}\label{T:2}

 Let $\shm\in Mod_{coh}(\shd_X)$. Then,  $$SS(\text{R}\shh\text{om}_{{\shd}_X} (\shm , Db_{X|N}))\subset \tilde{f}^{-1}(C^1_f(\shm)).$$
\end{corollary}

\section{Realification}
Let $f: X\to Y$ be a smooth morphism of complex manifolds. Let us denote by $\overline{Y}$ the underlying topological space $Y$ endowed with the complex conjugate structure and, identifying $Y$ to the diagonal $\Delta\subset Y\times\overline{Y}$ by the canonical embedding $\delta$, let us regard $Y\times \overline{Y}$ as a complexification of $Y$. Let us denote by $j: X\to X\times \overline{Y}$ the associated morphism and $p: X\times\overline{Y}\to X$ the associated projections.
\begin{proposition}\label{P:23}
For any $F\in
D^b_{\mathbb{R}-c}(\mathbb{C}_Y)$ and any $\shm\in D^b_{coh}(\shd_X)$ we have a canonical isomorphism:
\begin{equation}\label{E:25}
Rp_*\text{R}\shh\text{om}_{\shd_{X\times\overline{Y}}} (\underline{p}^{-1}\shm ,  t\shh \text{om}((f\times Id_{\overline{Y}})^{-1}R\delta_*F,\sho_{X\times \overline{Y}})[dim Y])$$
$$\simeq \text{R}\shh\text{om}_{\shd_X} (\shm , t\shh \text{om}(f^{-1}F,\sho_X)).
\end{equation}
\end{proposition}
\begin{proof}
We have $p j=Id_X$ and $j$ is proper so we may use the adjunction formula (7.4) of Theorem 7.2 of \cite{KS3} with respect to $j$, $f^{-1}F$ and $\underline{p}^{-1}\shm$ to obtain a natural isomorphism
\begin{equation}\label{E:205}\text{R}\shh\text{om}_{\shd_{X\times\overline{Y}}} (\underline{p}^{-1}\shm , t\shh \text{om}(Rj_* f^{-1}F,\sho_{X\times \overline{Y}})[dim Y])
\end{equation}
$$\simeq Rj_* \text{R}\shh\text{om}_{\shd_X} (\shm , t\shh \text{om}(f^{-1}F,\sho_X)).$$
On the other hand, the equality $\delta f=(f\times Id_{\overline{Y}}) j$ and Proposition 2.5.11 of \cite{KS1} entail the existence of a natural isomorphism such that, for any $F\in D^b_{\mathbb{R}-c}(\mathbb{C}_Y)$,
$Rj_*f^{-1}F\simeq (f\times Id_{\overline{Y}})^{-1}R\delta_*F$ which, composed with (\ref{E:205}), gives (\ref{E:25}).
\end{proof}
\begin{corollary}\label{C:206}
In the situation of Proposition \ref{P:23},  we have
$$SS(\text{R}\shh\text{om}_{\shd_X} (\shm , t\shh \text{om}(f^{-1}F,\sho_X)))$$ $$\subset p_{\pi}p_d(SS(\text{R}\shh\text{om}_{\shd_{X\times\overline{Y}}} (\underline{p}^{-1}\shm , t\shh \text{om}((f\times Id_{\overline{Y}})^{-1}R\delta_*F,\sho_{X\times \overline{Y}})))).$$
\end{corollary}
\begin{proof}
 This is a consequence of the properness of $p$ on the support of $$\text{R}\shh\text{om}_{\shd_{X\times\overline{Y}}} (\underline{p}^{-1}\shm , t\shh \text{om}((f\times id_{\overline{Y}})^{-1}R\delta_*F,\sho_{X\times \overline{Y}})).$$
\end{proof}

\begin{proposition}\label{C:207}

For any $\shm\in D^b_{coh}(\shd_X)$, we have
$$p_{\pi} p_d^{-1}((\widetilde{f\times Id_{\overline{Y}}})^{-1}(C^1_{f\times Id_{\overline{Y}}}(\underline{p}^{-1}\shm)))=\tilde{f}^{-1}(C^1_f(\shm)).$$
\end{proposition}
\begin{proof}

Let $q: T^*(X\times\overline{Y})\to T^*X$ be the projection.
Since $C^1_{f\times Id_{\overline{Y}}}(\underline{p}^{-1}\shm))\subset T^*(X|Y)\times \overline{Y}$, by Remark \ref{R:200} we have
 $$p_{\pi} p_d^{-1}((\widetilde{f\times Id_{\overline{Y}}})^{-1}(C^1_{f\times Id_{\overline{Y}}}(\underline{p}^{-1}\shm)))=q(\widetilde{(f\times Id_{\overline{Y}}})^{-1}(C^1_{f\times Id_{\overline{Y}}}(\underline{p}^{-1}\shm))).$$
 By  Lemma \ref{L:152} we have  $q((\widetilde{f\times Id_{\overline{Y}}})^{-1}(C^1_{f\times Id_{\overline{Y}}}(\underline{p}^{-1}\shm)))=\tilde{f}^{-1}C^1_f(\shm)$ which ends the proof.
\end{proof}

\section{Cauchy-Kowalevskaia theorem and propagation}

\subsection{Cauchy-Kowalevskaia theorem}
Let $X\subset \mathbb{C}^{n}\times \mathbb{C}^d$ be an open neighbourhood of $0$ of the form $X=Z\times Y$, with $Z\subset \mathbb{C}^n$ and $Y\subset\mathbb{C}^d$. $f$ will denote the projection $X\to Y$.

We shall consider a differential operator $P\in \shd_{X|Y}$ of the form
 \begin{equation}\label{E:208}
 P(z,y,D_z,D_y)=D_{z_1}^m+\underset{0\leq j<m}{\sum}a_{\alpha j}(z,y)D_{z'}^{\alpha}D_{z_1}^j,
 |\alpha|\leq m-j,
 \end{equation}

\noindent where $\alpha=(\alpha_2,...,\alpha_{n})\in \mathbb{N}^{n-1}$,  $D_{z'}^{\alpha}=D_{z_2}^{\alpha_2}...D_{z_{n}}^{\alpha_{n}}$ and  where the coefficients $a_{\alpha j}(z,y)$ are holomorphic in a neighbourhood $\Omega_0$ of $0$ in $X$.

Let  $W$ be an open subset in $Z$ and let $h\in \mathbb{C}$. We note $H_h$ and $H$, respectively, the hyperplane in $Z$ given by  $z_1=h$ and $z_1=0$. Let $\delta$ be a real positive number. According to the definition of J.M. Bony and P. Schapira in \cite{BS}, $W$ is $\delta-H_h-flat$ if $W$ is convex and for  any two given points  $z\in W$ and $\tilde{z}\in H_h$, the condition
 $$ |z_1-h|\geq\delta|z_j-\tilde{z_j}|_{j=2,...,n}$$ implies that  $\tilde{z}$ is contained in $W\cap H_h$.

Let now $\Omega$ be an open subset of $X$. We shall say that $\Omega$ is $Y-\delta-H_h-flat$ if, for any $y\in Y$, $f^{-1}(y)\cap\Omega$ is $\delta-H_h-flat$ in $f^{-1}(y)$.

In this situation we have:

\begin{proposition}\label{P:5}(\textbf{Precised Cauchy Theorem})\,\,
 Let $U$ be an open subanalytic set of $Y$.
  Let P be an operator of the form
(\ref{E:208}).
Then there exist an open neighbourhood $\Omega_0$ of $0\in X$ and $\delta>0$ such that, for any $h$, for any open subset $\Omega\subset\Omega_0$ which is  $Y-\delta-H_h-flat$, for any $$f\in \Gamma(\Omega;\sho_X^{t-Z\times U})$$ and for any $$(g)=(g_j)_{j=0,..., m-1} \in \Gamma(\Omega\cap H_h\times Y;\sho_{H_h\times Y}^{t-H_h \times U})^m $$ the Cauchy Problem
\begin{equation}\label{E:4}
Pu= f, \gamma(u)=(g),
\end{equation}
where $\gamma(u)=(u|_{H_h\times Y},..., D_{z_1}^m u|_{H_h\times Y})$, admits a unique solution $$u\in\
\Gamma(\Omega;\sho_X^{t-Z\times U}).$$
\end{proposition}
\begin{proof}
We shall adapt and follow step by step the proof of Theorem 2.4.3 of \cite{B}.
If $0\in U\cup (Y\setminus \overline{U})$ the assertion is the classical precised Cauchy Theorem. Hence we assume that $0\in \partial U$.
Let us start by considering as $\Omega_0$ an open neighborhood of $0$ where all the $a_{\alpha_j}$ are holomorphic.
We need to introduce some notations:

 Let $\tilde{\Omega}\subset\Omega_0$ be a relatively compact open subset.  For $(z,y)=(z_1, z', y)\in \tilde{\Omega}$ we set $$(i)\,\, d^{\prime}_{\tilde{\Omega}}(z,y)=\inf \{|| z'-\tilde{z}'||; (z_1,\tilde{z}',y)\in \complement\tilde{\Omega}\}.$$
For  $f\in \Gamma(\Omega;\sho_X^{t- Z\times U})$, for any relatively compact open subset $\tilde{\Omega}\subset \Omega$ and any $k\in \mathbb{N}$, we set
$$(ii)\,\,|f|_{\tilde{\Omega},k}=\underset{\tilde{\Omega}\cap Z\times U}{\sup} (|f(w)|d(w, \complement(Z\times U))^k d^{\prime}_{\tilde{\Omega}}(w)).$$

Therefore, the Cauchy data $f$  and $(g)$ satisfy the following condition:
for any relatively compact open subset $\tilde{\Omega}\subset \Omega$ (respectively $\tilde{\Omega^{\prime}}\subset \Omega\cap H_h\times Y$),
 there exists an integer $k\geq0$ such that $|f|_{\tilde{\Omega},k}< +\infty$   (respectively, for each $j=0,..., m-1$, $|g_j|_{\tilde{\Omega^{\prime}},k}< +\infty$).

 Notice that the operator $P$ may obviously be written in the form
\begin{equation}\label{E:3}
 P(z,y,D_z,D_y)=D_{z_1}^m+\underset{0\leq j<m-1}{\sum}D_{z'}^{\alpha}D_{z_1}^j a'_{\alpha_j}(z,y), |\alpha|\leq m-j,
 \end{equation}
 where the coefficients $a'_{\alpha_j}(z,y)$ are holomorphic in $\Omega_0$. It is then clear that, for any $\alpha_j$, any  relatively compact open subset $\tilde{\Omega}\subset\Omega$ and $f\in\Gamma(\Omega;\sho_X^{t-{Z\times U}})$, there exists $C_{\alpha_j}\geq 0$ such that
 \begin{equation}\label{E:8}
  |a'_{\alpha_j}f|_{\tilde{\Omega},k}< C_{\alpha_j} |f|_{\tilde{\Omega},k}.
  \end{equation}

 Up to shrinking $\Omega_0$ and $\delta$, any compact subset  $K$ of $\Omega$ is contained in a finite union of relatively compact open sets of the form

  $\tilde{\Omega}=\tilde{\Omega}_1\times \tilde{\Omega}_2$, with relatively compact open subsets $\tilde{\Omega}_1\subset Z$, $\tilde{\Omega}_2\subset Y$ and $\tilde{\Omega}_1$ being $\delta-H_{h}-flat$.

 We may therefore assume that $K$ is contained  in such a single $\tilde{\Omega}$.

 \begin{lemma}\label{L:251}
 Assume that the Cauchy data of problem (\ref{E:4})
satisfies
\begin{align}
\label{A:4}
    |f|_{\tilde{\Omega},k}<+\infty,\,\, & |g_j|_{\tilde{\Omega}^{\prime},k}<+\infty, j\leq m-1
     \end{align}
then the solution $u$ satisfies $|u|_{\tilde{\Omega},k}<+\infty$.
 \end{lemma}
 \begin{proof}
We may assume $h=0$ and then set $H=H_{0}$.
 We shall construct the solution $u$ of (\ref{E:4}) by successive approximation, defining recursively a sequence $$u_{\nu}\in \Gamma(\Omega;\sho_X^{t- {Z\times U}})$$ satisfying $|u_{\nu}|_{\tilde{\Omega},k}<+\infty$, by setting $u_0=0$ and $u_{\nu+1}$ being the solution of the Cauchy Problem
 \begin{equation}\label{E:5}
  D_{z_1}^m  u_{\nu +1}= f+ \underset{0\leq j\leq m-1,\,|\alpha|\leq m-j}{\sum}D_{z'}^{\alpha}D_{z_1}^j a'_{\alpha_j}u_{\nu},
   \end{equation}
 $$D_{z_1}^j   u_{\nu +1}| _{H\times Y} = g_j, j=0,..., m-1.$$

 Note that $\Omega$ being $Y-\delta-H-flat$, for any $(z_1, z', y)\in\Omega\cap Z\times U$ and any $t\in [0,1]$, the point $(tz_1, z', y)$ belongs to $\Omega\cap Z\times U$, hence problem $(\ref{E:5})$ is solvable.

 Let us note $v_{\nu}=u_{\nu +1 } -u_{\nu}$. We get $u_{\nu}=\underset{i=0,...,\nu -1}{\sum} v_{i}$,

 \begin{equation}\label{E:6}
   D_{z_1}^m v_0= f,
    D_{z_1}^j   v_0| _{H\times Y} = g_j, j=0,..., m-1
 \end{equation}
 and
 \begin{equation}\label{E:7}
   D_{z_1}^m  v_{\nu +1}=  \underset{0\leq j\leq m-1,\,|\alpha|\leq m-j}{\sum}D_{z'}^{\alpha}D_{z_1}^j a'_{\alpha j}v_{\nu},\,\,
 D_{z_1}^j   v_{\nu +1}| _{H\times Y} =0.
\end{equation}


 Let us fix a given $\nu$,   assume that $v_{\nu}\in \Gamma(\Omega;\sho_X^{t-{Z\times U}})$ and satisfies $|v_{\nu}|_{\tilde{\Omega},k}<+\infty$.    We shall prove that there exists a constant $C>0$ depending only on $P$ and $\Omega_0$, such that
 \begin{equation}\label{E:250}|v_{\nu+1}|_{\tilde{\Omega},k}\leq C\delta |v_{\nu}|_{\tilde{\Omega},k}.
 \end{equation}





 In view of (\ref{E:8}), it is enough to prove the following lemma:

 \begin{lemma}\label{L:25}
 There exists a constant $C>0$, only depending on $P$ and $\Omega_0$ such that, given $$\omega\in\Gamma(\Omega;\sho_X^{t -{Z\times U}})$$ satisfying  $|\omega|_{\tilde{\Omega},k}<+\infty$, the solution $v$ of the Cauchy problem

 $$D_{z_1}^m v=   D_{z'}^{\alpha}D_{z_1}^j\omega,
 \,(with\, |\alpha|\leq m-j, j<m)$$
$$D_{z_1}^j   v| _{H\times Y} =0, j<m$$

satisfies
\begin{equation}\label{E:14}
 |v|_{\tilde{\Omega},k} \leq C\delta |\omega|_{\tilde{\Omega},k}.
 \end{equation}
 \end{lemma}

\begin{proof}

 Remark that $v$ is also the solution of the Cauchy Problem
 \begin{equation}\label{E:10}
 D_{z_1}^{m-j}  v=   D_{z'}^{\alpha} \omega,
 D_{z_1}^j   v| _{H\times Y} =0, j\leq m-j-1.
\end{equation}
We have, for such  $\tilde{\Omega}$, $|\omega|_{\tilde{\Omega},k}=\underset{(z,y)\in \tilde{\Omega}_1\times (\tilde{\Omega}_2\cap U) }{\sup}(|\omega(z,y)|d(y, \complement U)^k d^{\prime}_{\tilde{\Omega}}(z,y)),$ therefore, for any $(z,y)\in \tilde{\Omega}_1\times (\tilde{\Omega}_2\cap U)$,
$$|\omega(z,y)|\leq |\omega|_{\tilde{\Omega},k} d(y, \complement U)^{-k}d^{\prime}_{\tilde{\Omega}}(z,y)^{-1}.$$

Let $z=(z_1,z')\in\tilde{\Omega}_1$ and let $y\in \tilde{\Omega}_2\cap U$.
By Cauchy integral formula, we get
\begin{equation}\label{E:11}
 |D_{z'}^{\alpha} \omega(z,y)|\leq|\omega|_{\tilde{\Omega},k}\, e^{|\alpha|}(|\alpha|+1)!d(y, \complement U)^{-k} d^{\prime}_{\tilde{\Omega}}(z,y)^{-|\alpha|-1}.
\end{equation}

Furthermore, we may assume $\delta< 1$ and  the diameter of $\Omega_0<1/m$.
Hence
 \begin{equation}\label{E:12}
 |D_{z'}^{\alpha} \omega(z,y)|\leq|\omega|_{\tilde{\Omega},k}\, e^{|\alpha|}(|\alpha|+1)!d(y, \complement U)^{-k} d^{\prime}_{\tilde{\Omega}}(z,y)^{-m+j-1}.
\end{equation}
 Since $\tilde{\Omega}_1$ is $\delta-H-flat$, for $t\in [0, 1]$,  $p_t:=(tz_1, z', y)$ belongs to $\tilde{\Omega}\cap Z\times U$ and  $$ d^{\prime}_{\tilde{\Omega}}(p_t)\geq d^{\prime}_{\tilde{\Omega}}(z,y)  +|\frac{(1-t)z_1}{\delta}|$$
 hence,
 \begin{equation}\label{E:13}
 |\int_0^1D_{z'}^{\alpha} z_1dt|\leq |\omega|_{\tilde{\Omega},k}\, e^{|\alpha|}(|\alpha|+1)!d(y, \complement U)^{-k} \frac{\delta}{m-j} d^{\prime}_{\tilde{\Omega}}(z,y)^{-m+j}.
 \end{equation}
 Iterating $m-j$ times this integration, we get
 $$|v(z,y))|\leq |\omega|_{\tilde{\Omega},k} e^{|\alpha|}(|\alpha|+1)!d(y, \complement U)^{-k} \frac{\delta ^{m-j}}{(m-j)!}d'_{\tilde{\Omega}}(z,y)^{-1},$$ hence
 $$|v|_{\tilde{\Omega},k}\leq |\omega|_{\tilde{\Omega},k} e^{|\alpha|}(|\alpha|+1)!\frac{\delta ^{m-j}}{(m-j)!},$$
 and we can choose $C= e^{|\alpha|}(|\alpha|+1)!$.
\end{proof}
 Remark that $v_0=G+\underset{j=0,...,m-1}{\sum} g_j (z',y)z_1^j /j!$ where $G$ denotes the mth-primitive of $f$, vanishing up to the order $m$ on $H\times Y$. By assumption, $|f|_{\tilde{\Omega},k}<+\infty, \, |g_j|_{\tilde{\Omega}\cap H\times Y,k}<+ \infty$.  Hence   $|v_0|_{\tilde{\Omega},k}<+ \infty$. Noting $C_0=|v_0|_{\tilde{\Omega},k}$, arguing by induction thanks to  (\ref{E:8}) and (\ref{E:250}) we get
 \begin{equation}\label{E:15}
 |v_{\nu +1} |_{\tilde{\Omega},k}\leq C_0 (C\delta)^{\nu +1}.
  \end{equation}
  We may also assume that $\delta<1/C.$
  Therefore the series $\sum v_{\nu}$ defines a holomorphic function $u$ in $\Omega\cap Z\times U$ satisfying  the estimate
  $|u|_{\tilde{\Omega},k}\leq C_0 \sum (C\delta)^{\nu}<+\infty.$
   Clearly, $u$ satisfies $Pu=f$ and $u$ is unique.

This ends the proof of Lemma \ref{L:251}.
 \end{proof}

 Proposition \ref{P:5} follows imediately from Lemma \ref{L:251}.
 \end{proof}

 We shall now obtain in Theorem \ref{T:401} below a tempered variant of Cauchy-Kowalevskaia-Kashiwara theorem (cf \cite{K1}).

  In the following results we consider $X$ and $Y$ complex manifolds, we assume $f:X\to Y$ be a given  smooth morphism and $\shm\in D^b_{coh}(\shd_X)$ be given. Let us note $V=X\times_Y T^*Y$.
 Let be given a submanifold  $H$ of $X$ non 1-characteristic for $\shm$ with respect to $f$ and let $f_H$ denote $f|_H$.

 \begin{corollary}\label{C:20}

 Assume $F=\mathbb{C}_U$ where $U$ is a  finite union of Stein open subanalytic relatively compact subsets of $Y$.

  Then the natural morphism in $D^{b}(\mathbb{C}_H)$
 \begin{equation}\label{E:2091}
 \text{R}\shh\text{om}_{{\shd}_X} (\shm ,t \shh \text{om} (f^{-1} F, \sho_X))|_H\to \text{R}\shh\text{om}_{{\shd}_H} (\shm_H ,t \shh \text{om} (f_H^{-1}F, \sho_H))
 \end{equation}
  is an isomorphism.
 \end{corollary}
 \begin{proof} The statement follows easily by the classical methods first reducing to the case of a single coherent module $\shm$ and then to the case where $H$ is an hypersurface and where  $\shm$ is defined by an operator in the conditions of  Proposition \ref{P:5}.   To end the proof one argues by induction on the number of open subanalytic sets whose union gives $U$, using the fact that  the family of Stein open subanalytic sets is stable under finite intersection.
\end{proof}
 \begin{corollary}\label{T:40}
 Assume that $Y$ is the complexification of a real analytic manifold $N$ and that $F\in D^b_{\mathbb{ R}-c}(\mathbb{C}_N).$
 Then the natural morphism in $D^{b}(\mathbb{C}_H)$
 \begin{equation}\label{E:209}
 \text{R}\shh\text{om}_{{\shd}_X} (\shm ,t \shh \text{om} (f^{-1} F, \sho_X))|_H\to \text{R}\shh\text{om}_{{\shd}_H} (\shm_H ,t \shh \text{om} (f_H^{-1}F, \sho_H))
 \end{equation}
  is an isomorphism.
\end{corollary}
 \begin{proof}

 We may assume the following situation in a neighbohood of $x$:
 \noindent $X=Z\times Y$, where  $Z$ is an open neigbourhood of $0\in \mathbb{C}^{n}$, $Y$ is an open neighbourhood of $0\in\mathbb{C}^d$ and $f$ is the projection.





 It is enough to prove the result for $F=\mathbb{C}_U$, with $U$ being an open subanalytic relatively compact subset of $N$.
Moreover we may assume that $N$ is $\mathbb{R}^d$ and $Y=\mathbb{C}^d$.
By the subanalytic Grauert's theorem (\cite{BMF}), we can write $U$ as  $U= \Omega\cap \mathbb{R}^d$,
where  $\Omega$ is Stein open subanalytic in $\mathbb{C}^d$. On the other hand, $\mathbb{C}^d\setminus \mathbb{R}^d$ is a finite union of Stein open subsets. The result then follows by Corollary \ref{C:20}.
\end{proof}

Keeping the same notations and assumptions, the preceding Corollaries and Proposition \ref {P:23} together entail:

\begin{theorem}\label{T:401}(Cauchy-Kowalevskaia-Kashiwara theorem)

 Let $F\in D^b_{\mathbb{ R}-c}(\mathbb{C}_X)$ and assume that $SS(F)\subset V$. Then the natural morphism in $D^{b}(\mathbb{C}_H)$
 \begin{equation}\label{E:209}
 \text{R}\shh\text{om}_{{\shd}_X} (\shm ,t \shh \text{om} ( F, \sho_X))|_H\to \text{R}\shh\text{om}_{{\shd}_H} (\shm_H ,t \shh \text{om} (F|_H, \sho_H))
 \end{equation}
  is an isomorphism.
\end{theorem}

\subsection{Propagation}

 Notice that, for  an open  subset  $\Omega$ of a complex manifold $X$,  $U$  a subanalytic open subset of $X$  and $p\in\partial \Omega$,  given  $u\in\Gamma(\Omega;\sho_X^{t-U})$,  the precise meaning of the assertion  that $u$ extends as a section of $\sho_X^{t-U}$ in  a neighbourhood of $p$ is that   there exist an open neighborhood  $\Omega'$ of $p$
  and a section $u'\in\Gamma(\Omega';\sho_X^{t-U})$, such that, as  holomorphic functions, $u'|_{\Omega'\cap\Omega\cap U}=u|_{\Omega'\cap\Omega\cap U}$.

 Let us now come back to the situation of Proposition \ref{P:5}. Hence  $X\subset \mathbb{C}^{n}\times \mathbb{C}^d$ is an open neighbourhood of $0$ of the form $X=Z\times Y$, for some open subsets $Z\subset \mathbb{C}^n$, $Y\subset\mathbb{C}^d$, $f$  denotes the projection $X\to Y$ and  $P\in \shd_{X|Y}$ satisfies (\ref{E:208}).
The following Lemma is an adaptation of Zerner's Lemma(\cite{Z}):
 \begin{lemma}\label{L:2}
Let $U$ be an open subanalytic set of $Y$. Let $\phi$ be a $C^{\infty}$ function in a neighbourhood of $0\in X$  such that $\phi(0)=0$ and $d\phi(0)=dz_1$. Let $\Omega=\{(z,y)\in X: \phi(z,y)<0\}.$ Assume that $Pu$ extends as a section of $\sho_X^{t-{Z\times U}}$ in  a neighbourhood of $0$.Then $u$ extends  (as a section
of $\sho_X^{t-{Z\times U}}$) in a neighbourhood of $0$.
 \end{lemma}
 \begin{proof}

In view of Proposition \ref{P:5} the proof is similar to that of Lemma 2.7 in \cite{TMF}. We may assume that $\phi$ is defined in $\Omega_0\subset \Omega$, with $\Omega_0$ satisfying the assumptions of Proposition \ref{P:5} with respect to a given $\delta$, and that $Pu$ extends to $\Omega_0$.  Then there exists $0<\epsilon\ll1$ and $R>0$ such that the open polydisc centered in $(-\epsilon, 0,0)$ and radius max $(R, \delta R)$ is contained in $\Omega_0$.
 Then, $$W_{\epsilon}=\{(z_1, z',y): |z_1+\epsilon|<\delta(R-||z'||), ||z'||<R, ||y||<R\}\subset\Omega_0$$ is $Y-\delta-H_{-\epsilon}-flat$ and is a neighbourhood of $0$. Again by Proposition  \ref{P:5} the solution of $Pu_{\epsilon}=Pu, \gamma(u_{\epsilon})=\gamma(u)$  defines a section of   $\sho_X^{t-{Z\times U}}$ on $W_{\epsilon}$.
 \end{proof}
 The preceding Lemma has a global version adapting Lemma 3.1.5 of \cite{S}:
 \begin{lemma}\label{L:20}
 Let $\omega$ and $\Omega$ be two convex  subsets in $X$, with $\Omega$ open, $\omega$ locally closed and $\omega\subset\Omega$. Assume that any real hyperplane whose conormal belongs to the closure of $\{\xi; \exists x\in \Omega, \tilde{f} (x,\xi)\in C^1_f(\shd_X/\shd_X P)\}$ which intersects $\Omega$ also intersects $\omega$. Then, if $u\in \Gamma(\omega; \sho_X^{t-{Z\times U}})$ is such that $Pu$ extends to $\Omega$ (as a section in  $\Gamma(\Omega; \sho_X^{t-{Z\times U}})$, then $u$ extends to $\Omega$.
 \end{lemma}

 \begin{remark}\label{L:3}

 Let $Z$, $Y$ be open sets respectively in $\C^n$ and $\C^d$  and let $f:Z\times Y\to Y$ be  the projection. Let $U$ be a subanalytic open set in $Y$.
Given  an open subset $\Omega$ of $X=Z\times Y$ such that, for any $y\in U$, $\Omega\cap f^{-1}(y)$ is  connected, and given   an open subset $\omega\subset \Omega$ such that, for any $y\in U$,  $\omega\cap f^{-1}(U)$ is non empty,

\textit{ if $u\in \Gamma(\Omega;\sho_X^{t-Z\times U)})$ vanishes on $\omega$, it follows that $u=0$.}

\noindent Indeed, regarding $u$ as a holomorphic function on $\Omega\cap f^{-1}(U)$, we may apply the analytic continuation principle to the restriction of $u$ to each $\Omega\cap f^{-1}(y)$, $y\in U$, as a holomorphic function in the $z$-variable.
Since $u|_{\omega\cap f^{-1}(y)}$ vanishes, $u|_{\Omega\cap f^{-1}(y)}$ vanishes and the statement follows.

Note that in the case of $f=Id:\C^d\to\C^d$, this result is trivial since the conditions above just mean that $U\subset\omega$ and that $u$, regarded as a holomorphic fonction on $\omega\cap U=U=\Omega\cap U$, vanishes.
 \end{remark}

 \begin{lemma}\label{L:4}(Propagation)
 Let $X$ and $Y$ complex manifolds and let $f:X\to Y$  be a smooth morphism. Let $U$ be a finite union of Stein subanalytic relatively compact open subsets of $Y$. Then, for any $\shm\in D^b_{coh}(\shd_X)$, we have
 \begin{equation}\label{E:17}
 SS(\text{R}\shh\text{om}_{{\shd}_X} (\shm , t \shh \text{om} (\mathbb{C}_{f^{-1} U}  , \sho_X)))\subset \tilde{f}^{-1}(C^1_f(\shm)).
\end{equation}
 \end{lemma}
 \begin{proof}
By induction on the number of open subsets whose reunion gives $U$, we may assume that $U$ is  a relatively compact connected Stein open set. Let $p=(x, \xi)\notin \tilde{f}^{-1}(C^1_f(\shm))$.
 As above, we may assume that
 $ x=0\in \mathbb{C}^{n+d}, X=Z\times Y$, where  $Z$ is an open neigborhood of $0\in \mathbb{C}^{n}$, $Y$ is an open neighbourhood of $0\in\mathbb{C}^d$, $f$ is the projection and $p= dz_1$. Moreover, we may assume that $\shm$ is of the form $\shd_X/ \shd_X P$ with
 \begin{equation}\label{E:19}
 P(z,y,D_z,D_y)=D_{z_1}^m+\underset{0\leq j<m}{\sum}a_{\alpha}(z,y)D_{z'}^{\alpha}D_{z_1}^j,
 |\alpha|\leq m
 \end{equation}
 where $\alpha=(\alpha_2,...,\alpha_{n})\in \mathbb{N}^{n-1}$,  $D_{z'}^{\alpha}=D_{z_2}^{\alpha_2}...D_{z_{n}}^{\alpha_{n}}$ and  where the coefficients $a_{\alpha}(z,y)$ are holomorphic in a neighbourhood $\Omega_0$ of $0$.
 Let $\shh$ denote the complex  $\text{R}\shh\text{om}_{{\shd}_X} (\shm , t \shh \text{om} (\mathbb{C}_{Z\times U}  , \sho_X))$.
 Let $\Omega$ and $ \delta$ be given by Proposition \ref{P:5}.
 We shall apply Proposition \ref{P:165}  and prove that, for $0<\epsilon\ll1$ and $0<R\ll1$, denoting by $V$ the open polydisc of center $(0,...,0)$ and radius  $R$, $B_R(0)$, by $F_{\varepsilon}=\{(z,y)\in X;\text{Re}
z_1>-\varepsilon\},$ by $L_{\epsilon}=\{(z,y); \text{Re}z_1=-\epsilon\}$ and by $\gamma$ the proper closed convex cone of $\mathbb{C}^{n+d}$
 $$
 \gamma=\{(z_1,z',y), \text{Re} z_1\leq-\delta( ||(z',y)||+|\text{Im}z_1|)\},$$  we have
 \begin{center}
$H^j(X; \mathbb{C}_{((z,y)+\gamma)\cap F_{\epsilon}} \otimes
\shh)=0, \text{for any}\, j\in \mathbb{Z}, \text{and any} (z,y)\in V.$
\end{center}

Remark that $p\in \text{int}\gamma^{oa}$.
 We may assume $\epsilon$ and $R$ small enough such that for $(z,y)\in V$, $((z,y)+\gamma)\cap\overline{F_{\epsilon}}\subset\Omega_0$.
 It is enough to prove that P defines an isomorphism in
 $$\shh^j(X; \mathbb{C}_{((z,y)+\gamma)\cap F_{\epsilon}} \otimes t \shh \text{om} (\mathbb{C}_{Z\times U}  , \sho_X) ).$$

We claim that
\begin{equation}\label{E:211}
\text{for}\, j\neq 1, \shh^j(X; \mathbb{C}_{((z,y)+\gamma)\cap F_{\epsilon}} \otimes t \shh \text{om} (\mathbb{C}_{Z\times U}  , \sho_X))=0.
\end{equation}
Set $F'_{\epsilon}=\{(z,y): Re z_1<-\epsilon\}$.

Case $j=0$:

 Let us prove that:
$$H^0 (X; \mathbb{C}_{((z,y)+\gamma)\cap F_{\epsilon}} \otimes t \shh \text{om} (\mathbb{C}_{Z\times U}  , \sho_X))=0.$$
  Let $s\in\Gamma(X; \mathbb{C}_{((z,y)+\gamma)\cap F_{\epsilon}} \otimes  t \shh \text{om} (\mathbb{C}_{Z\times U}  , \sho_X))$. Then there exists $(\tilde{z},\tilde{y})\in V$ such that $s$ extends as a section $s'$ of $\Gamma ((
\tilde{z},\tilde{y})+\text{int}(\gamma);
  t \shh \text{om} (\mathbb{C}_{Z\times U}  , \sho_X))$ with support in $\overline{F_{\epsilon}}$. Take $\Omega=(
\tilde{z},\tilde{y})+\text{int}(\gamma)$ and
$\omega=\Omega\cap F'_{\epsilon}$. Hence  $s$ vanishes on $\omega$. Since $\Omega$ is convex, for any $y\in U$, $f^{-1}(y)\cap \Omega$ is convex. Moreover  $f^{-1}(y)\cap \omega\neq \emptyset$ provided $\delta$ is small enough. 
By Remark \ref{L:3}  it follows that $s=0$.

Case $j>1$:

Set  $\overline{{F'}_{\epsilon}}=\{(z,y)\in X;\text{Re}
z_1\leq-\varepsilon\}$.
Since $(z,y)+\gamma,\, ((z,y)+\gamma)\cap \overline{F'}_{\epsilon}$ are closed convex, by Proposition \ref{P:300} we get, for $j\geq1$,
$$H^j(X; \mathbb{C}_{((z,y)+\gamma)} \otimes t \shh \text{om} (\mathbb{C}_{Z\times U}  , \sho_X))=0,$$  $$H^j (X; \mathbb{C}_{((z,y)+\gamma)\cap \overline{{F'}_{\epsilon}}} \otimes t \shh \text{om} (\mathbb{C}_{Z\times U}  , \sho_X))=0.$$
Then (\ref{E:211}) follows by the distinguished triangle
\begin{center}$R\Gamma(X;  \mathbb{C}_{((z,y)+\gamma))\cap F_{\epsilon}} \otimes t \shh \text{om} (\mathbb{C}_{Z\times U}  , \sho_X))\to$
$R\Gamma(X; \mathbb{C}_{(z,y)+\gamma}\otimes t\mathcal{H}\text{om} (\mathbb{C}_{Z\times U}  , \sho_X))\to R\Gamma(X;\mathbb{C}_{((z,y)+\gamma)\cap \overline{{F'}_{\epsilon}}} \otimes t \shh \text{om} (\mathbb{C}_{Z\times U}, \sho_X))\underset {+1}{\to}$\end{center}



 Therefore, it is enough to prove that $P$ defines an isomorphism in $$\shh^1(X; \mathbb{C}_{((z,y)+\gamma)\cap F_{\epsilon}} \otimes t \shh \text{om} (\mathbb{C}_{Z\times U}  , \sho_X)).$$

 It is sufficient to prove that $P$ defines an isomorphism on $$\shh^1((z,y)+\gamma)\cap \overline{F_{\epsilon'}}, t \shh \text{om} (\mathbb{C}_{Z\times U}  , \sho_X))$$ for $\epsilon>\epsilon'>0$. By successive application of Proposition \ref{P:300} and Remark \ref{L:3}, we get
 \begin{center}$\shh^1((z,y)+\gamma)\cap \overline{F_{\epsilon'}}, t \shh \text{om} (\mathbb{C}_{Z\times U}  , \sho_X))\simeq\frac{\Gamma(((z,y)+\gamma)\cap L_{\epsilon'}; t \shh \text{om} (\mathbb{C}_{Z\times U}  , \sho_X))}{\Gamma(((z,y)+\gamma)\cap \overline{F_{\epsilon'}}; t \shh \text{om} (\mathbb{C}_{Z\times U}  , \sho_X))}$\end{center}


 Note that any given $(z,y)\in L_{\epsilon'}$ admits a fundamental system of $Y-\delta-H_{-\epsilon'}$-flat neighborhoods.
 So, given $v\in \Gamma(((z,y)+\gamma)\cap L_{\epsilon'}; t\shh\text{om}(\mathbb{C}_{Z\times U}  ,\sho_X))$, by Proposition \ref{P:5} we solve the equation $Pu=v$ in a neighborhood of $H_{-\epsilon'}$ and then we extend it to a neighborhood of $((z,y)+\gamma)\cap\overline{{F'}_{\epsilon}}$ by Lemma \ref{L:20}.
 This proves the surjectivity of $P$.
 To prove the injectivity we follow \cite{TMF}, which is based on $(5.1.4)$ and $(5.1.5)$ of \cite{KS1}.
 For each $a=(z_0,y_0)\in V\cap(\overline{F_{\epsilon'}}\setminus  L_{\epsilon'})$ one constructs a family of open subsets $\{\Omega_t(a)\}_{t\in\R^+}$, such that:
\begin{enumerate}
\item{$\Omega_t(a)\subset a+\mathrm{int}\gamma$,}
\item{$\Omega_t(a)\cap L_{\epsilon'}=(a+\mathrm{int}\gamma)\cap L_{\epsilon'}$,}
\item{$\Omega_t(a)=\underset{r<t}{\cup}(a)$,}
\item{$\partial\Omega_t(a)$ is smooth real analytic,}
\item{$Z_t(a):=(\underset{s>t}{\cap}\overline{\Omega_s(a)\setminus\Omega_t(a)})\cap\overline{F_{\epsilon'}}\subset \partial \Omega_t(a)$ and the conormal to $\Omega_t(a)$ at the points in $Z_t(a)$ is non 1-characteristic for $P$,}
\item{$(\underset{t>0}{\cup}\Omega_t(a))\cap\overline{F_{\epsilon'}}\subset(a+\mathrm{int}\gamma)\cap\overline{F_{\epsilon'}}$,}
\item{$(\underset{t>0}{\cap}\Omega_t(a))\cap\overline{F_{\epsilon'}}\subset(a+\mathrm{int}\gamma)\cap{L_{\epsilon'}}$.}
\end{enumerate}
Remark that, for $v=(1,0,...,0; 0,...,0)$  the family $$\{(z,y)+\rho v+\mathrm{int}\gamma\cap\overline{F_{\epsilon'}}\}_{\rho>0}$$ forms a neighborhood system of $(z,y)+\gamma \cap\overline{F_{\epsilon'}}$, and the family $$\{\Omega_t((z,y)+\rho v)\cap \overline{F_{\epsilon'}}\}_{\rho>0, t>0}$$ forms a neighborhood system of $((z,y)+\gamma)\cap {L_{\epsilon'}}$. Let $f\in \Gamma(((z,y)+\gamma)\cap L_{\epsilon'};  t\shh\text{om}(\mathbb{C}_{Z\times U}  ,\sho_X))$ such that $Pf=g$ extends to a neighborhood of $((z,y)+\gamma)\cap \overline{F_{\epsilon'}}.$ Choose $\rho$ and $t_0$ such that $f$ is defined in $\Omega_{t_{0}}((z,y)+\rho v)\cap \overline{F_{\epsilon'}}$ and such that $g$ is defined in $((z,y)+\rho v +\mathrm{int}\gamma)\cap \overline{F_{\epsilon'}}.$ By the assumptions and Lemma \ref{L:2} we can extend $f$ to $\Omega_{t'}((z,y)+\rho v)\cap \overline{F_{\epsilon'}}$ for some $t'>t_0$ and this procedure gives an extension of $g$ to $((z,y)+\gamma)\cap \overline{F_{\epsilon'}}$.
 \end{proof}


\begin {proposition}\label{P:10}
 Let $X$ and $Y$ complex manifolds and let $f:X\to Y$  be a smooth morphism. Assume that $Y$ is a complexification of a real manifold $N$. Let $F$ belong to $D^b_{\mathbb{R}-c}(\mathbb{C}_N)$. Then (\ref{E:1}) holds for $G=f^{-1}F$.
 \end{proposition}
\begin{proof}

It is enough to prove the result for $F=\mathbb{C}_U$ for an open subanalytic relatively compact subset of $N$.
Moreover we may assume that $N$ is $\mathbb{R}^d$ and $Y=\mathbb{C}^d$.
By the subanalytic Grauert's theorem (cf.\cite{BMF}), we can write $U$ as  $U=\Omega\cap \mathbb{R}^d$,
where $\Omega$ is Stein open subanalytic in $\mathbb{C}^d$. By Lemma \ref{L:4}, (\ref{E:1}) holds for $  \Omega$. On the other hand, $\mathbb{C}^d\setminus \mathbb{R}^d$ is a finite union of Stein open subsets, hence (\ref{E:1}) holds for $\Omega\setminus \mathbb{R}^d$ which ends the proof.
\end{proof}


\begin{proof}[Proof of Theorem \ref{T:1}]

The problem being of local nature on $X$, we may assume that $G\simeq  f^{-1}F$ for some $F\in D^b_{\mathbb{R}-c}(\mathbb{C}_Y).$
By Corollary \ref{C:206} we have
\begin{center}$SS(\text{R}\shh\text{om}_{\shd_X} (\shm , t\shh \text{om}(f^{-1}F,\sho_X))) \subset p_{\pi}p_d(SS(\text{R}\shh\text{om}_{\shd_{X\times\overline{Y}}} (\underline{p}^{-1}\shm , t\shh \text{om}(f\times id_{\overline{Y}})^{-1}R\delta_*F,\sho_{X\times \overline{Y}}))).$\end{center}
where $p:X\times\overline{Y}\to X$ denotes the projection. The result is then an immediate consequence of Proposition \ref{P:10} together with Corollary \ref{C:207}.
\end{proof}

\begin{proof}[Proof of Corollary \ref{T:2}]
It is an immediate consequence of Proposition \ref{P:10}.
\end{proof}


In the case of a product, $f$ being a projection, there is a particular class of $\shd_X$-modules for which the estimate in Theorem \ref{T:1} is easily improved.

\begin{corollary}\label{C:2}
Assume that $X=Z\times Y$, let  $p: X\to Z$, $f:X\to Y$  and $q:T^*Z\times T^*Y\to T^*Z$ denote the  projections. Let $\shm\in D^b_{coh}(\shd_Z)$ and let $G\in D^b_{\R-c}(\C_X)$ be given and satisfying $SS(G)\subset V$. Then
$$SS(\text{R}\shh\text{om}_{\shd_X} (\underline{p}^{-1}\shm , t\shh \text{om}(G,\sho_X)))\subset q^{-1}(Char(\shm)).$$
\end{corollary}

\begin{proof}
 Reducing first to $\shm$ in degree zero, and secondly to $\shm\simeq\shd_Z/\mathcal{J}$, for a coherent  ideal of $\shd_Z$, one easily checks the equality

\begin{center}$\tilde{f}^{-1}(C^1_f(\underline {p}^{-1}\shm))=q^{-1}(Char(\shm)).$\end{center}
 \end{proof}

\appendix

\section{Microlocalization}

\subsection{1-microcharacteristic varieties}
Let $V$ be a smooth involutive subamanifold of $T^*X\setminus T^*_X X$, let $\she_X(m)_{m\in\Z}$ denote the filtration of $\she_X$ by the order,  and let $\mathcal{J}_V$ denote the subsheaf of $\she_X$ of microdifferential operators of order at most $1$ whose symbol  of order $1$ vanishes on $V$. Let $I_V\subset \sho_{T^*X}$ be the defining ideal of $V$.  One denotes $\she_V$ the sub-sheaf of rings of $\she_X$ generated by $\mathcal{J}_V$,
$$\she_V=\cup_{m\geq0} \mathcal{J}^m.$$
For $P\in \she_V\cap \she_X(m)$, $P\notin \she_X(m-1)$, one defines $\sigma_V^1(P)$ as being the image of $\sigma(P)$ in $I_V^m/I_V^{m+1}.$ Hence $\sigma_V^1(P)$ defines a holomorphic function on $T_V(T^*X).$

Let $\shm$ be a coherent $\she_X$-module. By the constructions of \cite{TMF1} and \cite{L},  one associates to $\shm$ the 1-microcharacteristic variety along $V$, a conic analytic subset of $T_V(T^*X)$, satisfying, for $\shm=\she_X/\mathcal{J}$, $\mathcal{J}$ a coherent ideal of $\she_X$ ,
$$C^1_V(\shm)=\{p\in T_V(T^*X); \sigma^1_V(P)(p)=0, \forall P\in \mathcal{J}\cap \she_V\}.$$
The definition extends naturally to $D^b_{coh}(\she_X).$

Recall that the canonical $1$-form on $T^*X$  defines a section  $T^*X\underset{H}{\hookrightarrow} T^*(T^*X)$. We also note $\rho_V$ the composition of the Hamiltonian isomorphism $T^*(T^*X)\to T(T^*X)$ with the canonical projection $$V\underset {T^*X}{\times} T(T^*X)\to T_V(T^*X).$$

Let now $f:X\to Y$ be a smooth morphism and set $V=T^*Y\times_Y X$. Set $V^{\cdot}=V\setminus T_X^*X$. Then $\pi^{-1}\shd_{X|Y}\subset \she_{V^{\cdot}}$, more precisely, it is equal to $\pi^{-1}\shd_X\cap \she_{V^{\cdot}}$.

Let now $\shm\in D^b_{coh}(\shd_X)$ and set $\tilde{\shm}=\she_X\otimes_{\pi^{-1}\shd_X} \pi^{-1}\shm.$

\begin{proposition}\label{L:302}

 $H^{-1}\rho_V^{-1}(C^1_{V^{\cdot}}(\tilde{\shm}))=\tilde{f}^{-1}(C^1_f(\shm)).$
\end{proposition}
 \begin{proof}
We may assume that $\shm$ is concentrated in degree zero and that  $\shm=\shd_X/\mathcal{L}$ for a coherent ideal $\mathcal{L}\subset \shd_X$. Moreover, since the problem is local, we may assume that $f:X=Z\times Y\to Y$ is the projection. Taking local coordinates $(z,y)$ in $X$ and the associated canonical coordinates $(z,y; \xi,\eta)$ on $T^*X$, we get $V=\{(z,y;\xi,\eta)\in T^*X, \xi=0\}.$ The result then follows because,  for $P\in \shd_{X|Y}$, one has
$$\sigma^1_{V^{\cdot}}(P)(z,y,0,\eta;\tilde{\xi})=\sigma^1_{X|Y}(P)(z,y;\tilde{\xi}).$$
\end{proof}

\subsection{Estimate for $t\mu hom$}
In the results below we adapt the techniques in \cite{KMS}.

Let us recall the construction
of the functor $t\mu hom(\cdot,{\sho}_X)$ of tempered
microlocalization of (\cite{A}) (also we refer \cite{SG} where it is proved that    $t\mu hom(F,{\sho}_X)$ is an object of $D(\she_X)$).

Let $\tilde{X}^{\mathbb{C}}$ be the complex normal deformation of $X\times{X}$ along
 the diagonal $\Delta$ which we
 identify with $X$ by the first
 projection $p_1$.  This gives an identification of
$TX$ with the normal bundle $T_{\Delta}(X\times{X})$. Let
$t:\tilde{X}^{\mathbb{C}}\to{\mathbb{C}}$ and
$p:\tilde{X}^{\mathbb{C}}\to{X\times{X}}$ be  the
 canonical maps,
let $\tilde{\Omega}$ be $t^{-1}({\mathbb{C}}-\{0\})$ and
$\Omega=t^{-1}(\mathbb{R}^{+})\subset\tilde{\Omega}$.
 Let
$p_2:X\times{X}\to X$ be the  second projection.

Consider the following diagram of morphisms:
\begin{equation}\label{E:7}
TX\simeq{T_{\Delta}(X\times{X})}\overset{i}{\hookrightarrow}\tilde{X}^{\mathbb{C}}\overset{j}{\hookleftarrow}\Omega=t^{-1}(\mathbb{R}^{+})
\end{equation}
Let $\tilde{p}:\tilde{\Omega}\to{X\times{X}}$,~
 be the
 restriction of $p$. Denote by $\overline{p}_{1}$ the composition $p_{1}\circ{p}$ and by
$\overline{p}_{2}$ the composition $p_{2}\circ{p}$.

Under these notations, $t\nu hom(F,{\sho}_X)$ is defined by
 $$t\nu hom(F,{\sho}_X)
={i^{-1}R{\shh}om_{{\shd}_{\tilde{X}^{\mathbb{C}}}}({\shd}_{\tilde{X}^{\mathbb{C}}\underset{\overline{p}_{1}}{\to}X}
,t\shh om(\overline{p}_{2}^{-1}F\otimes{\mathbb{C}_{\Omega}},{\sho}_{\tilde{X}^{\mathbb{C}}}))}.$$
Let $D^b_{\R^+}(\mathbb{C}_{TX})$
(resp. $D^b_{\R^+}(\mathbb{C}_{T^*X})$)
be the derived category
of complexes of sheaves on $TX$ (resp. $T^*X$)
with conic cohomologies.
We denote by the symbol $\,\,\widehat{}\,\,$ the
Fourier-Sato Transform from
$D^b_{\R^+} (\mathbb{C}_{TX})$ to
$D^b_{\R^+} (\mathbb{C}_{T^*X})$. Then,
by definition, $t\mu hom(F,{\sho}_X)=t\nu hom(F,{\sho}_X)^{\,\widehat{}}$ .
 Let us recall that under the identification of $T^*(TX)$ with
$T^*(T^*X)$ by  the Hamiltonian isomorphism
we have  $SS(F)=SS(F)^{\,\widehat{}}$ for any
$F\in{D^b_{\R^+}(\mathbb{C}_{TX})}$.

  Remark that for any coherent ${\shd}_X$-module
${\shm}$, one has
 \begin{align}\label{E:49}
 & R{\shh}om_{\pi^{-1}{\shd}_X}(\pi^{-1}{\shm},t \mu hom(F,{\sho}_X))\\
&\simeq
{ R{\shh}om_{\tau^{-1}{\shd}_X}(\tau^{-1}{\shm},t\nu h om(F,{\sho}_X))}^{\,\widehat{}}.\notag
\end{align}
As before, let $\tilde{\shm}=\she_X\otimes_{\pi^{-1}\shd_X} \pi^{-1}\shm$.
\begin{theorem}\label{T:300}
Let $\shm\in D^b_{coh}(\shd_X)$. Let $f:X\to Y$ be a smooth morphism and let $V=T^*Y\times _Y X$. Let $G\in D^b_{\R-c}(\C_X)$ and assume that $SS(G)\subset V$.
Then
 \begin{equation}\label{E:301}
 SS(R\shh om_{\she_X}(\tilde{\shm}, t\mu hom(G, \sho_X)))\subset \rho_V^{-1}(C^{1}_{V^{\cdot}}(\tilde{\shm})).
  \end{equation}

\end{theorem}
\begin{proof}

By (\ref{E:49}) it is sufficient to prove the analogue of (\ref{E:301}) with $t\mu
hom(G,{\sho}_X)$  replaced by $t\nu hom(G,{\sho}_X)$. Moreover, the assertion being local, we may assume that $f:X=Z\times Y\to Y$ is the projection.
We have
\begin{align}&R\shh om_{\tau^{-1}\shd_X}(\tau^{-1}\shm,t\nu hom(G,\sho_X))\notag\\
&\simeq
{i^{-1}R{\shh}om_{{\shd}_{\tilde{X}^{\C}}}(\underline{\overline{p}_1}(\shm),
t\shh om(\overline{p}_{2}^{-1}G\otimes{\mathbb{C}_{\Omega}},
\sho_{\tilde{X}^{\mathbb{C}}}))}\notag
\end{align}

 By Proposition 6.6.2 of \cite{KS1}, we have an
inclusion

\begin{center}\label{E:55}
$SS(~i^{-1}R{\shh}om_{{\shd}_{\tilde{X}}^{\mathbb{C}}}(\underline{\overline{p}_1}(\shm),
t \shh om(\overline{p}_{2}^{-1}G\otimes{\mathbb{C}_{\Omega}},{\sho}_{\tilde{X}^{\mathbb{C}}})))\subset{i^{\sharp}(SS(R{\shh}om_{{\shd}_{\tilde{X}^{\mathbb{C}}}}(\underline{\overline{p}_1}(\shm),
t\shh om(\overline{p}_{2}^{-1}G\otimes{\mathbb{C}_{\Omega}},{\sho}_{\tilde{X^{\mathbb{C}}}}))))}$.
\end{center}

Endow $X\times{X}$
with the system of local coordinates
$(x,x')$, so that
 $\Delta\subset{X\times{X}}$ is defined by $x=x'$. Under
the change of coordinates : $x=x'$, $y=x-x'$, $\Delta$ will be
defined by $y=0$.  The coordinates $(x,y)$ induce in $\tilde{X}^{\mathbb{C}}$  coordinates $(t,x,y)$,
such that
\begin{equation}\label{E:51}
\mbox{$p(t,x,y)
=(x,x-ty)$,  $\overline{p}_1(t,x,y)=x$, and
$\overline{p}_2(t,x,y)=x-ty$}
\end{equation}
We may assume that the local coordinate system $x$ is of the form $x=(z,u)$
such that $f\colon X\to Y$ is given by $(z,u)\mapsto u$.
Then $V=X\times_YT^*Y$ is given by
$V=\{(z,u;\xi,\eta);\xi=0\}$.
Set $y=(z^{\prime},u^{\prime})$.
Let $f_2$ denote the smooth morphism $\tilde{X}^{\C}\to Y\times \C$ given by $$f_2(t,x,y)=(f \circ  \overline{p}_2(t,x,y),t)=(u-tu^{\prime}, t).$$

Let $(t,z,u,z^{\prime},u^{\prime};\tau,\xi,\eta,\xi^{\prime},\eta^{\prime})$ be the associated coordinates of
$T^*(\tilde{X}^{\mathbb{C}})$.

Since $\overline{p}_2^{-1}G\otimes \C_{\Omega}$ is locally of the form $f_2^{-1}(F\boxtimes \C_{\{\Re t>0\}})$, by Theorem \ref{T:1} we have:
\begin{equation}\label{E:56}
SS(R{\shh}om_{{\shd}_{\tilde{X}^{\mathbb{C}}}}
(\underline{\overline{p}_1}(\shm),t\shh om(\overline{p}_{2}^{-1}G
\otimes{\mathbb{C}_{\Omega}},{\sho}_{\tilde{X}^{\mathbb{C}}})))
\subset \tilde{f_2}^{-1}(C^1_{f_2}(\underline{\overline{p}_1}(\shm))).
\end{equation}
 By  (\ref{E:55}) and (\ref{E:56}) we get
 \begin{equation}\label{E:57}
SS(R{\shh}om_{\tau^{-1}{\shd}_X}(\tau^{-1}{\shm},t\nu
hom(G,{\sho}_X))) \subset{i}^{\sharp}(\tilde{f_2}^{-1}(C^1_{f_2}(\underline{\overline{p}_1}(\shm)
))).
\end{equation}
Therefore it is enough to prove   the inclusion
\[{i}^{\sharp}(\tilde{f_2}^{-1}(C^1_{f_2}(\underline{\overline{p}_1}(\shm)
)))\subset \rho_V^{-1}(C^1_V(\tilde{\shm})).\]

We may assume that $\shm$ is concentrated in degree zero and is of the form $\shm=\shd_X/\mathcal{J}$ for a coherent ideal $\mathcal{J}$. In this case we have
$$\underline{\overline{p}_1}\shm\simeq \frac{\shd_{\tilde{X}^{\C}}}{\shd_{\tilde{X}^{\C}}\mathcal{J}+\langle D_{z^{\prime}}, D_{u^{\prime}},D_t\rangle}.$$
Since $\langle D_{z^{\prime}}, \mathcal{J}\cap\shd_{X|Y}\rangle\subset \shd_{\tilde{X}^{\C}|Y\times\C}$ we have \[\tilde{f_2}^{-1}(C^1_{f_2}(\underline{\overline{p}_1}(\shm)))\subset\{(t,z,u,z^{\prime},u^{\prime};\tau,\xi,\eta,\xi^{\prime},\eta^{\prime}); \xi^{\prime}=0, \sigma^1_{X|Y}(P)(z,u;\xi)=0,\] \[ \forall P\in\mathcal{J}\cap \shd_{X|Y}\}.\]
Hence
\[{i}^{\sharp}(\tilde{f_2}^{-1}(C^1_{f_2}(\underline{\overline{p}_1}(\shm)
)))\subset\{(z,u,z^{\prime}, u^{\prime}, \xi,\eta, 0, \eta^{\prime})\in T^*(TX);\] \[
\exists (t_n, z_n, u_n, z^{\prime}_n, u^{\prime}_n; \tau_n,\xi_n, \eta_n, 0, \eta^{\prime}_n),\]
\[ t_n\to 0, z_n\to z, u_n\to u, \xi_n\to\xi, \eta_n\to \eta, \eta^{\prime}_n\to\eta^{\prime}, |t_n\tau_n|\to 0, \] \[\sigma^1_{X|Y}(P)(z_n, u_n; \xi_n)=0, \forall P\in\mathcal{J}\cap\shd_{X|Y}\}\]
\[\subset\{(z,u,z^{\prime}, u^{\prime}, \xi,\eta, 0, \eta^{\prime}); \sigma^1_{X|Y}(P)(z, u;\xi)=0, \forall P\in\mathcal{J}\cap\shd_{X|Y}\}.\]

Recall that the identification of $T^*(TX)$ to $T(T^*X)$ is described by:
$$T^*(TX)\ni(z,u, z^{\prime},u^{\prime}; \xi, \eta, \xi^{\prime}, \eta^{\prime})\leftrightarrow (z,u,\xi^{\prime},\eta^{\prime}; z^{\prime}, u^{\prime},\xi,\eta)\in T(T^*X).$$Hence
\[{i}^{\sharp}(\tilde{f_2}^{-1}(C^1_{f_2}(\underline{\overline{p}_1}(\shm)
)))\subset \{(z,u,0,\eta^{\prime}; z^{\prime}, u^{\prime},\xi,\eta)\in T(T^*X);\sigma^1_{X|Y}(P)(z,u;\xi)=0\}\] \[ =\{(z,u,0,\eta^{\prime}; z^{\prime}, u^{\prime},\xi,\eta)\in T(T^*X); \sigma^1_{V^{\cdot}}(P)(z,u,0,\eta^{\prime};\xi)=0, \forall P\in\mathcal{J}\cap\she_{V^{\cdot}}\}=\] \[ \rho_V^{-1}(C^1_{V^{\cdot}}(\shm)).\]
\end{proof}


\begin{thebibliography}{15}
\bibitem{A}
E. Andronikof, {\em Microlocalization temp\' er\' ee}, Suppl\'emment  au Bull.Soc.Math.France, m\' emoire 57,
\textbf{122}, 2, (1994).
\bibitem{BMF} D. Barlet and T. Monteiro Fernandes, {\em Grauert's theorem for subanalytic open sets in real analytic manifolds}, Studia Mathematica, \textbf{204},3, 265-274, (2011).
\bibitem{Be}
O. Berni, {\em A vanishing theorem for formal cohomology of perverse sheaves}, Journal of Funct. Anal., \textbf{158}, 2, 267-288, (1998).
\bibitem{B}
J.M. Bony, {\em Propagation des singularit\' es diff\'erentiables pour une classe d'op\' erateurs diff\'erentiels \`a coeffitients analytiques}, Ast\'erisque, \textbf{34-35 }, 43-92, (1976).
\bibitem{BS}
J.M. Bony and P. Schapira, {\em Propagation des singularit\' es analytiques pour les solutions des \'equations aux d\'eriv\'ees partielles}, Ann. Inst. Fourier, \textbf{26}, 81-140, (1976)

\bibitem{AAT}A. D'Agnolo and F.Tonin, {\em Cauchy Problem for hyperbolic $\shd$-modules with regular singularities}, Pacific J. Math. \textbf{184}, 1-22 (1998).
\bibitem{SG}
S. Guillermou, {\em DG-Methods for Microlocalization}, Publ. RIMS, Kyoto Univ., \textbf{1}, 99-140, (2011).


\bibitem{G}
H. Grauert, {\em On Levi's Problem and the embedding of real analytic manifolds}, Annals of Maths, \textbf{68}, 2, 460-472 (1958).
\bibitem{Ho}
L. H\"ormander, {\em An introduction to Complex Analysis in several variables}, North-Holland, \textbf{7}, (1973).
\bibitem{Ho1}
L. H\"ormander, {\em The analysis of linear partial differential
operators II,} Grundlehren der Math. Wiss.\, \textbf{257}, Springer
Verlag, (1983).



\bibitem{KT}H. Koshimizu and K. Takeuchi, {\em Extension theorems for distribution solutions to $\shd$-modules with regular singularities}, Proceedings of the AMS, \textbf{128}, 6, 1685-1690 (2000).
\bibitem{K1}
M. Kashiwara, {\em D-Modules and Microlocal Calculus},  Translations of Mathematical Monographs, \textbf{ 217}, AMS (2003).
\bibitem{K2}
 M. Kashiwara, {\em The Riemann-Hilbert problem for holonomic systems}, Publ. RIMS, Kyoto Univ, \textbf{20}, 219-315 (1984).




\bibitem{KS1}
M. Kashiwara and P. Schapira, {\em Sheaves on manifolds,}
Grundlehren der Math. Wiss., \textbf{292}, Springer Verlag (1990).

\bibitem{KS3}
M. Kashiwara and P. Schapira, {\em Moderate and formal cohomology associated with constructible sheaves,} Bull.Soc. Math. France, M\' emoire  \textbf{94} (1996).
\bibitem{KMS}
M. Kashiwara, T. Monteiro Fernandes and P. Schapira, {\em  Micro-support and Cauchy problem for temperate solutions of regular D-modules,} Portugaliae Mathematica, Vol 58, pp. 485-504 (2001).
\bibitem{KO}
M. Kashiwara and T. Oshima, {\em Systems of differential equations with regular singularities and their boundary values problem}, Annals of Math., \textbf{106}, 145-200 (1977).

\bibitem{L} Y. Laurent, {\em Th\'eorie de la deuxi\`eme microlocalization}, Progr. Math., Birkhauser,\textbf{53} (1985).


\bibitem{TMF}
T. Monteiro Fernandes,  {\em Propagation of the irregularity of a microdifferential system}, Publ. RIMS, Kyoto Univ., \textbf{37}, 2, 119-139, (2001)

\bibitem{TMF1}
T. Monteiro Fernandes {\em Probl\`eme de Cauchy pour les syst\`emes microdiff\'erentiels}, Ast\'erisque, \textbf{140-141}, 135-220, (1986).


\bibitem{S-K-K} M. Sato, T. Kawai and M.
Kashiwara, {\em Hyperfunctions and pseudodifferential equations,}
Lecture Notes in Math., Springer \textbf{287}, 265-529 (1973).


\bibitem{S}
P. Schapira, {\em Microdifferential systems in the complex domain},
Grundlehren der Math. Wiss., Vol. \textbf{269}, Springer-Verlag,
(1985).

\bibitem{Z}
M. Zerner, {\em Domaine d'holomorphie des fonctions v\'erifiant une \'equation aux d\'eriv\'ees partielles}, C. R. Acad. Sci. Paris S\'erie A, \textbf{272}, 1646-1648 (1971).
\end{thebibliography}
\end{document}